\documentclass[12pt]{article}
\usepackage{amssymb,amsmath}
\usepackage{mathrsfs}
\usepackage{graphicx}

\newtheorem{Thm}{Theorem}
\newcommand{\arch}{\mathrm{arcosh}}

\begin{document}

\begin{center}
\large{\bf SPREAD OF A CATALYTIC BRANCHING RANDOM WALK ON A MULTIDIMENSIONAL LATTICE}
\end{center}
\vskip0,5cm
\begin{center}
Ekaterina Vl. Bulinskaya\footnote{ \emph{Email address:} {\tt
bulinskaya@yandex.ru}}$^,$\footnote{The work is partially supported
by Dmitry Zimin Foundation ``Dynasty'' and RFBR grant 14-01-00318.}
\vskip0,2cm \emph{Lomonosov Moscow State University}
\end{center}
\vskip1cm

\begin{abstract}

For a supercritical catalytic branching random walk on $\mathbb{Z}^d$, $d\in\mathbb{N}$, with an arbitrary finite catalysts set we study the spread of particles population as time grows to infinity. Namely, we divide by $t$ the position coordinates of each particle existing at time $t$ and then let $t$ tend to infinity. It is shown that in the limit there are a.s. no particles outside the closed convex surface in $\mathbb{R}^d$ which we call the propagation front and, under condition of infinite number of visits of the catalysts set, a.s. there exist particles on the propagation front. We also demonstrate that the propagation front is asymptotically densely populated and derive its alternative representation. Recent strong limit theorems for total and local particles numbers established by the author play an essential role. The results obtained develop ones by Ph.Carmona and Y.Hu (2014) devoted to the spread of catalytic branching random walk on $\mathbb{Z}$.

\vskip0,5cm {\it Keywords and phrases}: branching random walk,
supercritical regime, spread of population, propagation front,
many-to-one lemma.

\vskip0,5cm 2010 {\it AMS classification}: 60J80, 60F15.

\end{abstract}

\section{Introduction}
Theory of branching processes is a vast and rapidly developing area of probability
theory having a multitude of applications (see, e.g., monographs \cite{HJV_05} and \cite{KA_15}). A branching
process is intended to describe evolution of population of individuals (particles) which could be genes,
bacteria, humans, clients waiting in a queue etc. A special section of that theory is constituted by processes
in which particles besides producing offspring also move in space. Such a scenario where the motion of a
particle is governed by random walk is named a branching random walk (for random walk, see, e.g., books
\cite{Lawler_Limic_10} and \cite{Borovkov_book_13}). One of the most natural and intriguing questions related to
branching random walk is how the particles population spreads in the space whenever it survives. Within the last
decades a lot of attention has been paid to that question in the framework of different models of branching
random walk on integer lattices or in Euclidean space. One can list publications since the paper \cite{Biggins_78}
till  numerous recent works, for instance, papers \cite{Bertacchi_Zucca_15},
\cite{Comets_Popov_07}, \cite{Le-Gall_Lin_15}, \cite{Mallein_15} and the monograph \cite{Shi_LNM_15}. However,
those results only slightly concern the model of {\it catalytic branching random walk} (CBRW) on
$\mathbb{Z}^d$, $d\in\mathbb{N}$, with a finite set of catalysts, which is considered here. A specific trait
of CBRW is its non-homogeneity in space, i.e. particles may produce offspring only at selected ``catalytic''
points of $\mathbb{Z}^d$ and the set of these points where catalysts are located is finite. This model is
closely related to the so-called parabolic Anderson problem (see, e.g., \cite{GKM_07}) and requires special
research methods.

Study of different variants of CBRW goes back to more than 10 years
(see, e.g., \cite{AB_00} and \cite{VTY}), although most of papers in
this research domain have been published recently, see, for
instance, \cite{Y_TPA_11}, \cite{HTV_12},
\cite{Molchanov_Yarovaya_2012}, \cite{B_TrudyMIAN_13}, \cite{DR_13},
\cite{B_JTP_14} and \cite{Carmona_Hu_2014}. A lot of them analyze
asymptotic behavior of total and local particles numbers as time
tends to infinity and only few investigate the spread of CBRW.
Analysis of the mean total and local particles numbers implemented
in the most general form in \cite{B_TPA_15} as well as the strong
and weak limit theorems established in \cite{B_Doklady_15} shows
that CBRW can be classified as supercritical, critical and
subcritical like ordinary branching processes and only in the
supercritical regime the total and local particles numbers grow
jointly to infinity. For this reason, it is of primary interest to
consider spread of particles population in supercritical CBRW.

The following advances in the study of CBRW spread have been
achieved. The paper \cite{Carmona_Hu_2014} devoted to CBRW on
$\mathbb{Z}$ reveals that the maximum of CBRW (i.e. the rightmost
particle location) increases asymptotically linearly in time tending to infinity.
Its authors employ the many-to-few lemma proved in general form
in \cite{Harris_Roberts_16}, martingale technique and renewal theorems. A similar assertion for catalytic
branching Brownian motion on $\mathbb{R}$ with binary fission and a single catalyst is established in
\cite{Bocharov_Harris_14} among other results. S.A.Molchanov and E.B.Yarovaya in their papers such as
\cite{Molchanov_Yarovaya_2012} study the spread of CBRW with binary fission and symmetric random walk on
$\mathbb{Z}^d$ by employing the operator theory methods for symmetric evolution operator. Note that in
\cite{CKMV_09}
the authors apply the continuous-space counterpart of such CBRW to modeling of homopolymers.

The main aim of our paper is to study the spread of CBRW on
$\mathbb{Z}^d$ for arbitrary positive integer $d$. In contrast to
the one-dimensional case where the maximum of CBRW on $\mathbb{Z}$
was investigated, one cannot directly extend the same approach to
multidimensional lattices and employ the fundamental martingale
techniques as in \cite{Carmona_Hu_2014}.  The point is that the
concept of maximum is indefinite for CBRW on $\mathbb{Z}^d, d>1$.
Were the random walk symmetric and catalysts positioned symmetrically, as well as the starting point of CBRW
be at the origin, then it would be sufficient to consider the maximum of the norm of particle locations or
the maximal displacement of a particle, similar to \cite{Mallein_15}. However, in more general setting it is of
interest to understand not only how far a particle can move from the origin but also in which direction such
displacement takes place.
So, in this paper, we introduce the concept of the propagation front
$\mathcal{P}\subset\mathbb{R}^d$ of the particles population as
follows. Divide by $t$ the position coordinates of each particle
existing in CBRW at time $t$ and let $t$ tend to infinity. Then in
the limit there are a.s. no particles outside the set bounded by the
closed surface $\mathcal{P}$ and, under condition of infinite number
of visits of catalysts, a.s. there exist particles on $\mathcal{P}$.
Thus, under this condition, non-random set $\mathcal{P}$ asymptotically separates the
a.s. population areal and its a.s. void environment. Moreover, we
establish that each point of $\mathcal{P}$ is a limiting point for
the normalized particles positions in CBRW and derive an alternative
representation for the propagation front $\mathcal{P}$. The latter
formula allows us to evaluate directly (without any computer
simulation) the set $\mathcal{P}$ for a number of examples in the
end of the paper. The proofs involve many-to-one formula, renewal theorems for systems of
renewal equations, martingale change of measure, convex analysis,
large deviation theory and the coupling method. We also essentially
base on recent investigation in \cite{B_TPA_15} of the mean total
and local particles numbers in CBRW as well as on the strong and
weak limit theorems for those quantities in \cite{B_Doklady_15}.

The paper is organized as follows. In
section~\ref{s:notation_results} we recall the necessary background
material and formulate three new theorems. Theorem~\ref{T:spread}
establishes the asymptotically linear pattern of the population
propagation  with respect to time growing to infinity.
Theorem~\ref{T:one_point} demonstrates that the set $\mathcal{P}$ is
asymptotically densely populated. Theorem~\ref{T:front} provides an
alternative representation for the front $\mathcal{P}$. In
section~\ref{s:spread} we establish both Theorems~\ref{T:spread} and
\ref{T:one_point} casting the proof into 5 steps.
Section~\ref{s:front} is devoted to the proof of
Theorem~\ref{T:front} and consideration of five examples. The first
example is related to CBRW on $\mathbb{Z}$ and we derive a result of
\cite{Carmona_Hu_2014} as a special case. Examples 2a, 2b and 2c
illustrate the spread of CBRW on $\mathbb{Z}^2$ in cases of
nearest-neighbor random walk, non-symmetric random walk and
non-symmetric random walk with unbounded jump sizes. Example 3
illustrates the spread of CBRW on $\mathbb{Z}^3$.

\section{Notation, main results and discussion}\label{s:notation_results}
Let us recall the description of CBRW on
$\mathbb{Z}^d$. At the initial time ${t=0}$ there is a
single particle that moves on $\mathbb{Z}^d$
according to a continuous-time Markov chain
${\bf S}=\{{\bf S}(t),t\geq0\}$ generated by the infinitesimal matrix
${Q=(q({\bf x},{\bf y}))_{{\bf x},{\bf y}\in\mathbb{Z}^d}}$. When this particle hits a finite set of
catalysts $W=\{{\bf w}_1,\ldots,{\bf w}_N\}\subset\mathbb{Z}^d$, say at the site ${\bf w}_k$,
it spends there random time having the exponential distribution with
parameter $\beta_k>0$. Afterwards the particle either branches or
leaves the site ${\bf w}_k$ with probabilities $\alpha_k$ and $1-\alpha_k$
($0\leq\alpha_k<1$), respectively. If the particle branches (at the
site ${\bf w}_k$), it dies and just before the death produces a random
non-negative integer number $\xi_{k}$ of offsprings located at the
same site ${\bf w}_k$. If the particle leaves ${\bf w}_k$, it jumps to the site
${\bf y}\neq {\bf w}_k$ with probability $-(1-\alpha_k)q({\bf w}_k,{\bf y})q({\bf w}_k,{\bf w}_k)^{-1}$
and continues its motion governed by the Markov chain ${\bf S}$. All
newly born particles are supposed to behave as independent copies of
their parent.

We assume that the Markov chain ${\bf S}$ is irreducible and the
matrix $Q$ is conservative (i.e. $\sum\nolimits_{{\bf
y}\in\mathbb{Z}^d}{q({\bf x},{\bf y})}=0$ where $q({\bf x},{\bf
y})\geq0$ for ${\bf x}\neq{\bf y}$ and $q({\bf x},{\bf
x})\in(-\infty,0)$ for any ${\bf x}\in\mathbb{Z}^d$). We employ the
standard assumption of existence of a finite derivative $f_k'(1)$,
that is the finiteness of $m_k:={\sf E}{\xi_k}$, for any
$k=1,\ldots,N$. Let $\mu(t)$ be the total number of  particles
existing in CBRW at time $t\geq0$ and the local particles numbers
$\mu(t;{\bf y})$ be the quantities of particles located at separate
points ${\bf y}\in\mathbb{Z}^d$ at time $t$.

While in \cite{Carmona_Hu_2014} the authors considered a
discrete-time CBRW we are interested in continuous-time process
since in the latter case we are able to employ directly new results
of \cite{B_TPA_15} and \cite{B_Doklady_15}. It is worthwhile to note
that in discrete-time and continuous-time settings most of
asymptotic results turn out to be the same modulo constants.
Moreover, in contrast to \cite{Carmona_Hu_2014} in this paper we
consider a variant of CBRW where there is an additional parameter
$\alpha_k$ governing the proportion between ``branching'' and
``walking'' of a particle located at each catalyst point ${\bf
w}_k$. However, as shown, e.g., in \cite{Y_TPA_11} and
\cite{B_JTP_14}, introducing of additional parameters does not
influence the asymptotic results for CBRW accurate up to constants.
At last, whereas in \cite{Carmona_Hu_2014} the underlying random
walk on $\mathbb{Z}$ is constructed as a cumulative sum of i.i.d.
random variables, in a similar way we assume that the underlying
random walk (i.e. our CBRW without branching) is space-homogeneous.
Due to the mentioned additional parameters it means that (see, e.g.,
\cite{Y_TPA_11})
\begin{equation}\label{condition1}
q({\bf x},{\bf y})=q({\bf x}-{\bf y},{\bf 0})=q({\bf 0},{\bf y}-{\bf x})\quad\mbox{and}\quad\beta_k=q/(1-\alpha_k),
\end{equation}
for ${\bf x},{\bf y}\in\mathbb{Z}^d$ and $k=1,\ldots,N$, where $q:=-q({\bf 0},{\bf 0})\in(0,\infty)$. Thus, our
investigation here can be considered as a development of the study of spread of CBRW initiated in
\cite{Carmona_Hu_2014} for one-dimensional case.

To formulate the main results of the paper let us introduce
additional notation. As usual, let all random elements be defined on
the same probability space $(\Omega,\mathcal{F},{\sf P})$.   The
index ${\bf x}$ in expressions of the form ${\sf E}_{\bf x}$ and
${\sf P}_{\bf x}$ marks the starting point of either CBRW or the
random walk ${\bf S}$ depending on the context. We temporarily
forget that there are catalysts at some points of $\mathbb{Z}^d$ and
consider only the motion of a particle on $\mathbb{Z}^d$ in
accordance with Markov chain ${\bf S}$ with generator $Q$ and
starting point ${\bf x}$. The conditions imposed on elements $q({\bf
x},{\bf y})$, ${\bf x},{\bf y}\in\mathbb{Z}^d$, allow us to use an
explicit construction of the random walk on $\mathbb{Z}^d$ with
generator $Q$ (see, e.g., Theorem~1.2 in \cite{Bremaud_99}, Ch.~9,
Sec.~1). According to this construction ${\bf S}$ is a regular jump
process with right continuous trajectories and, for transition times
of the process $\tau^{(0)}:=0$ and
$\tau^{(n)}:=\inf\left\{t\geq\tau^{(n-1)}:{\bf S}(t)\neq {\bf
S}(\tau^{(n-1)})\right\}$, $n\ge1$, the following statement holds.
Random variables
$\left\{\tau^{(n+1)}-\tau^{(n)}\right\}_{n=0}^{\infty}$ are
independent and each of them has exponential distribution with
parameter $q$. Denote by $P=\{P(t),t\geq0\}$ the Poisson process
constructed by means of the random sequence
$\{\tau^{(n+1)}-\tau^{(n)}\}_{n=0}^{\infty}$, i.e. $P$ is the
Poisson process with intensity $q$. Let ${\bf Y}^i$ be the value of
the $i$th jump of the random walk ${\bf S}$ ($i=1,2,\ldots$). In
view of that Theorem~1.2 in \cite{Bremaud_99}, Ch.~9, Sec.~1, the
random variables ${\bf Y}^1,{\bf Y}^2,\ldots$ are i.i.d., have
distribution ${\sf P}({\bf Y}^1={\bf y})=q({\bf 0},{\bf y})/q$,
${\bf y}\in\mathbb{Z}^d$, ${\bf y}\neq{\bf 0}$, and do not depend on
the sequence $\{\tau^{(n+1)}-\tau^{(n)}\}_{n=0}^{\infty}$. In other
words, the formula
\begin{equation}\label{S(t)=representation}
{\bf S}(t)={\bf x}+\sum_{i=1}^{P(t)}{\bf Y}^i
\end{equation}
holds true (as usual, $\sum_{i\in\varnothing}{\bf Y}^i=0$) where ${\bf x}$ is the initial state of the
Markov chain ${\bf S}$. Due to this equality it is not difficult to show that ${\bf S}$ is a process with
independent increments. In what follows we consider the version of
the process ${\bf S}$ constructed in such a way.

Set
$$\tau_{\bf x}:=\mathbb{I}({\bf S}(0)={\bf x})\inf\{t\geq0:{\bf S}(t)\neq{\bf x}\},$$
i.e. the stopping time $\tau_{\bf x}$ (with respect to the natural filtration $(\mathcal{F}_t,t\geq0)$ of
process ${\bf S}$) is the time of the first exit from the starting point ${\bf x}$ of the random walk. As usual,
$\mathbb{I}(A)$ stands for the indicator of a set $A\in\mathcal{F}$. Clearly,
${\sf P}_{\bf x}(\tau_{\bf x}\leq t)=1-e^{-qt}$, ${\bf x}\in\mathbb{Z}^d$, $t\geq0$. Let
$$_T\overline{\tau}_{{\bf x},{\bf y}}:=\mathbb{I}({\bf S}(0)={\bf x})\inf\{t\geq0:{\bf S}(t+\tau_{\bf x})={\bf y},{\bf S}(u)\notin T,\tau_{\bf x}\leq u< t+\tau_{\bf x}\}$$
be the time elapsed from the exit moment of this Markov chain (in
other terms, particle) out of starting state ${\bf x}$ till the
moment of first hitting point ${\bf y}$ whenever the particle
trajectory does not pass the set $T\subset\mathbb{Z}^d$. Otherwise,
we put ${_T\overline{\tau}_{{\bf x},{\bf y}}=\infty}$. Extended
random variable $_T\overline{\tau}_{{\bf x},{\bf y}}$ is called
\emph{hitting time} of state ${\bf y}$ \emph{under taboo} on set $T$
after exit out of starting state ${\bf x}$ (see, e.g.,
\cite{B_SPL_14}). Denote by $_T\overline{F}_{{\bf x},{\bf y}}(t)$,
$t\geq0$, the improper cumulative distribution function of this
extended random variable and let $_T\overline{F}_{{\bf x},{\bf
y}}(\infty):=\lim_{t\to\infty}{_T\overline{F}_{{\bf x},{\bf
y}}(t)}$. Whenever the taboo set $T$ is empty, expressions
$_{\varnothing}\overline{\tau}_{{\bf x},{\bf y}}$ and
$_{\varnothing}\overline{F}_{{\bf x},{\bf y}}$ are shortened as
$\overline{\tau}_{{\bf x},{\bf y}}$ and $\overline{F}_{{\bf x},{\bf
y}}$. Mainly we will be interested in the situation when $T=W_k$
where $W_k:=W\setminus\{{\bf w}_k\}$, $k=1,\ldots,N$.

Further
$$F^{\ast}(\lambda):=\int\nolimits_{0-}^{\infty}{e^{-\lambda t}\,d{F(t)}},\quad\lambda\geq0,$$ denotes
the Laplace transform of a cumulative distribution function $F(t)$,
$t\geq0$, with support located on non-negative semi-axis. For
$j=1,\ldots,N$ and $t\geq0$ set $G_j(t):=1-e^{-\beta_j t}$. In
\cite{B_TPA_15} there was introduced a matrix function $D(\lambda)$
with values in irreducible matrices of size $N\times N$ for each
$\lambda\geq0$. Namely, $D(\lambda)=(d_{i,j}(\lambda))_{i,j=1}^N$
where
$$d_{i,j}(\lambda)=\delta_{i,j}\alpha_im_iG^{\ast}_i(\lambda)+(1-\alpha_i)G^{\ast}_i(\lambda)
{_{W_j}\overline{F}^{\ast}_{{\bf w}_i,{\bf w}_j}(\lambda)}$$ and $\delta_{i,j}$
is the Kronecker delta. According to Definition~$1$ in
\cite{B_TPA_15} CBRW is called {\it supercritical} if the Perron root (i.e.
positive eigenvalue being the spectral radius) $\rho(D(0))$ of the
matrix $D(0)$ is greater than $1$. Then in view of monotonicity of
all elements of matrix function $D(\cdot)$ there exists the solution
$\nu>0$ of equation $\rho(D(\lambda))=1$. As Theorem $1$ in
\cite{B_TPA_15} shows, just this positive number $\nu$ specifies the
rate of exponential growth of the mean total and local particles
numbers (in the literature devoted to population dynamics and
branching processes one traditionally speaks of Malthusian
parameter). More precisely, ${\sf E}_{\bf x}\mu(t)\sim A({\bf x})e^{\nu t}$ and
${\sf E}_{\bf x}\mu(t;{\bf y})\sim a({\bf x},{\bf y})e^{\nu t}$ as $t\to\infty$ (the explicit
formulae for functions $A(\cdot)$ and $a(\cdot,\cdot)$ are given in
\cite{B_TPA_15}). Exactly these means play the role of normalizing factors in Theorems~3 and 4 of \cite{B_Doklady_15} devoted to the strong and weak convergence of vectors of the total and local particles numbers in supercritical CBRW as time
grows to infinity. In the given paper we concentrate on just a supercritical CBRW on $\mathbb{Z}^d$.

Let $N(t)\subset\mathbb{Z}^d$ be the (random) set of particles existing in CBRW at time $t\geq0$.
For a particle $v\in N(t)$, denote by ${\bf X}_v(t)$ its position at time $t$. Introduce the set of infinite
number of visits of catalysts by
$$\mathcal{I}=\left\{\omega: \limsup_{t\to\infty}\{v\in N(t):{\bf X}_v(t)\in W\}\neq\varnothing\right\}\in\mathcal{F}.$$
The behavior of CBRW on this set complement $\mathcal{I}^c$ is trivial since for $t\geq t_0(\omega)$ large enough
either CBRW dies out or CBRW constitutes the system of some random walks starting respectively from
${\bf X}_v(\omega,t_0)$, $v\in N(t_0)$, at time $t_0$. The supercritical regime of CBRW guarantees that
${\sf P}(\mathcal{I})>0$ (see, e.g., Theorem~4 of \cite{B_Doklady_15}).

Assume that the function
\begin{eqnarray}\label{H(s)_definition}
H({\bf s})&:=&\sum_{{\bf x}\in\mathbb{Z}^d}e^{\langle {\bf s},{\bf x}\rangle}q({\bf 0},{\bf x})=\sum_{{\bf x}\in\mathbb{Z}^d}\left(e^{\langle {\bf s},{\bf x}\rangle}-1\right)q({\bf 0},{\bf x})\\
&=&q\left(\sum_{{\bf x}\in\mathbb{Z}^d,\,{\bf x}\neq{\bf 0}}e^{\langle {\bf s},{\bf x}\rangle}\frac{q({\bf 0},{\bf x})}{q}-1\right)=q\left({\sf E}e^{\langle{\bf s},{\bf Y}^1\rangle}-1\right)\nonumber
\end{eqnarray}
is finite for any ${\bf s}\in\mathbb{R}^d$ where $\langle\cdot,\cdot\rangle$ stands for the inner product of vectors. This assumption is Cram\'{e}r's condition for the jump value ${\bf Y}^1$ satisfied in $\mathbb{R}^d$. It is easy to check that the Hessian of $H$ is positive definite and, consequently, $H$ is a convex function. Put also $\mathcal{R}=\left\{{\bf r}\in\mathbb{R}^d:H({\bf r})=\nu\right\}$.

At last, let
\begin{equation}\label{definition_O}
\mathcal{O}_{\varepsilon}:=\{{\bf x}\in\mathbb{R}^d:\langle {\bf x},{\bf r}\rangle>\nu+\varepsilon\;\;\mbox{for at least one}\;\;{\bf r}\in\mathcal{R}\},\quad\varepsilon>0,
\end{equation}
\begin{equation}\label{definition_Q}
\mathcal{Q}_{\varepsilon}:=\{{\bf x}\in\mathbb{R}^d:\langle{\bf x},{\bf r}\rangle<\nu-\varepsilon\;\mbox{for any}\;{\bf r}\in\mathcal{R}\},\quad\varepsilon\in(0,\nu),
\end{equation}
$\mathcal{O}:=\mathcal{O}_0$, $\mathcal{Q}:=\mathcal{Q}_0$ and $\mathcal{P}:=\partial\mathcal{Q}=\partial\mathcal{O}$ where $\partial\mathcal{S}$ stands for the border of a set $\mathcal{S}\subset\mathbb{R}^d$. Note that each set $\mathcal{Q}_{\varepsilon}$, $\mathcal{Q}$ or $\mathcal{P}\cup\mathcal{Q}$ is convex as an intersection of half-spaces (see, e.g., Theorem~2.1 of \cite{Rockafellar_70}).

\begin{Thm}\label{T:spread}
Let conditions \eqref{condition1} and \eqref{H(s)_definition} be satisfied for supercritical CBRW on $\mathbb{Z}^d$. Then, for any ${\bf x}\in\mathbb{Z}^d$ and  $t\to\infty$, we have
\begin{equation}\label{T:assertion_1}
{\sf P}_{\bf x}\left(\omega:\forall\varepsilon>0\;\exists t_0=t_0(\omega,\varepsilon)\;\mbox{s.t.}\;\forall t\geq t_0\;\mbox{and}\;\forall v\in N(t),\;{\bf X}_v(t)/t\notin\mathcal{O}_{\varepsilon}\right)=1,
\end{equation}
\begin{equation}\label{T:assertion_2}
{\sf P}_{\bf x}\left(\left.\omega:\!\forall\varepsilon\in(0,\nu)\exists t_1=t_1(\omega,\varepsilon)\;\mbox{s.t.}\;\forall t\geq t_1\;\exists v\in N(t),\;{\bf X}_v(t)/t\notin\mathcal{Q}_{\varepsilon}\right|\mathcal{I}\right)\!=\!1
\end{equation}
where the sets $\mathcal{O}_{\varepsilon}$ and $\mathcal{Q}_{\varepsilon}$ are defined in formulae \eqref{definition_O} and \eqref{definition_Q}, respectively.
\end{Thm}

Theorem~\ref{T:spread} means that if we divide the position
coordinates of each particle existing in CBRW at time $t$ by $t$ and
then let $t$ tend to infinity, then in the limit there are a.s. no
particles outside the set $\mathcal{P}\cup\mathcal{Q}$ and under
condition of infinite number of visits of catalysts there are a.s.
particles on $\mathcal{P}$. In this sense it is natural to call the
border $\mathcal{P}$ the {\it propagation front} of the particles
population. The following theorem refines assertion
\eqref{T:assertion_2} of Theorem~\ref{T:spread} and states that each
point of $\mathcal{P}$ can be considered as a limiting point for the
normalized particles positions in CBRW.

\begin{Thm}\label{T:one_point}
Let conditions of Theorem~\ref{T:spread} be satisfied. Then, for each ${\bf y}\in\mathcal{P}$, one has
$${\sf P}_{\bf x}\left(\left.\omega:\forall t\geq0\;\exists v_{\bf y}=v_{\bf y}(t,\omega)\in N(t)\;\mbox{such that}\;\lim_{t\to\infty}\frac{{\bf X}_{v_{\bf y}}(t)}{t}={\bf y}\right|\mathcal{I}\right)=1.$$
\end{Thm}

It follows from the definition of set $\mathcal{P}$ that
$$\mathcal{P}=\{{{\bf x}\in\mathbb{R}^d:{\langle{\bf x},{\bf r}\rangle\!\leq\!\nu}}\;\mbox{for all}\; {{\bf r}\in\mathcal{R}}\;\mbox{and}\;{\langle{\bf x},{\bf r}\rangle\!=\!\nu}\;\mbox{for at least one}\;{{\bf r}\in\mathcal{R}}\}.$$
Theorem~\ref{T:front} yields one more way to find the propagation front $\mathcal{P}$.
\begin{Thm}\label{T:front}
The set $\mathcal{P}$ can be also specified as $\mathcal{P}=\{{\bf
z}({\bf r}):{\bf r}\in\mathcal{R}\}$ where
$${\bf z}({\bf r})=\frac{\nu}{\langle\nabla H({\bf r}),{\bf r}\rangle}\nabla H({\bf r}).$$
\end{Thm}

This theorem allows us to evaluate directly (without any computer simulation) set $\mathcal{P}$ for a number of examples in section \ref{s:front} of the paper. Moreover, it follows from the proof of Theorem~\ref{T:front} that the definition of $\mathcal{P}$ can be refined as
$$\mathcal{P}=\{{{\bf x}\in\mathbb{R}^d:{\langle{\bf x},{\bf r}\rangle=\nu}}\;\mbox{for a single}\;{{\bf r}\in\mathcal{R}}\;\mbox{and}\;{\langle{\bf x},{\bf r}\rangle<\nu}\;\mbox{for other}\;{{\bf r}\in\mathcal{R}}\}.$$

Note that our new results show that the particles population spreads asymptotically linearly on $\mathbb{Z}^d$ with respect to growing time and the form of the propagation front does not depend on the number of catalysts and their locations but depends only on the value of the Malthusian parameter $\nu$ and the function $H(\cdot)$ characterizing the random walk. In other words, in our limit theorems the normalizing factor of the particles positions is equal to $t^{-1}$ and does not depend on the dimension of the lattice.

Remark that in \cite{Molchanov_Yarovaya_2012} there is also used the concept of the propagation front of CBRW with binary fission and symmetric random walk on $\mathbb{Z}^d$, namely, $\Gamma_t=\left\{{\bf y}={\bf y}(t)\in\mathbb{Z}^d:{\sf E}_{\bf 0}\mu(t;{\bf y})<C\right\}$ where $C$ is some positive constant. Moreover, in the framework of our terminology there is shown that $\Gamma_t=t\left(\mathcal{P}\cup\mathcal{O}\right)$. On the other hand, our formula \eqref{P(A^c_t)estimate} and its counterpart in case of multiple catalysts imply also that ${\sf E}_{\bf 0}\mu(t;{\bf y})<C$ for some positive constant $C>0$ and any ${\bf y}={\bf y}(t)\in t\left(\mathcal{P}\cup\mathcal{O}\right)$. Note also that we concentrate on almost sure results and impose less restrictions on the model than other researchers.

The present study became feasible due to the many-to-few formulae derived in the most general form in \cite{Harris_Roberts_16} and then applied to CBRW with a single catalyst in \cite{DR_13}. In a similar way one can obtain the following many-to-one formula for CBRW with several catalysts
\begin{equation}\label{many-to-one_formula}
{\sf E}_{\bf x}\sum_{v\in N(t)}g({\bf X}_v(t))={\sf E}_{\bf x}g({\bf S}(t))\prod_{k=1}^N\exp\{\alpha_k \beta_k(m_k-1)L(t;{\bf w}_k)\}
\end{equation}
where ${\bf x}\in\mathbb{Z}^d$, $L(t;{\bf y}):=\int_0^t\mathbb{I}({\bf S}(u)={\bf y})\,du$, ${\bf y}\in\mathbb{Z}^d$, $t\geq0$, is the local time of the random walk ${\bf S}$ at level ${\bf y}$ and $g:\mathbb{R}^d\rightarrow\mathbb{R}$ is a measurable function.
As noted above, in this paper we also employ renewal theorems for systems of renewal equations, martingale change of measure, convex analysis, large deviation theory and the coupling method. We essentially use results in \cite{B_TPA_15} on the mean total and local particles numbers in CBRW as well as the strong and weak limit theorems for those quantities established in \cite{B_Doklady_15}.

\section{Proof of Theorems~\ref{T:spread} and \ref{T:one_point}}\label{s:spread}
In this section we establish both Theorems~\ref{T:spread} and \ref{T:one_point} devoted to the study of spread of CBRW on $\mathbb{Z}^d$. For the sake of clarity of exposition their common proof is divided into 5 steps.

\vskip0.2cm{\it Step 1.} At the first step we assume that $W=\{{\bf w}_1\}$ with ${\bf w}_1={\bf 0}$ and the starting point of CBRW is ${\bf 0}$ as well.
Let us derive the first statement \eqref{T:assertion_1} of Theorem \ref{T:spread} for this case.

Fix ${\bf r}\in\mathcal{R}$. Let $\varepsilon>0$ and put $\mathcal{O}_{{\bf r},\varepsilon}:=\{{\bf x}\in\mathbb{R}^d:\langle{\bf x},{\bf r}\rangle>\nu+\varepsilon\}$.
According to \eqref{many-to-one_formula} we have
$${\sf P}_{\bf 0}\left(\exists v\in N(t):\;{\bf X}_v(t)\in t\mathcal{O}_{{\bf r},\varepsilon}\right)={\sf P}_{\bf 0}\left(\sum\limits_{v\in N(t)}\mathbb{I}\{{\bf X}_v(t)\in t\mathcal{O}_{{\bf r},\varepsilon}\}\neq0\right)$$
$$\leq{\sf E}_{\bf 0}\sum\limits_{v\in N(t)}\mathbb{I}\{{\bf X}_v(t)\in t\mathcal{O}_{{\bf r},\varepsilon}\}$$
$$={\sf E}_{\bf 0}\left(\mathbb{I}\{{\bf S}(t)\in t\mathcal{O}_{{\bf r},\varepsilon}\}\exp\{\alpha_1\beta_1(m_1-1)L(t;{\bf 0})\}\right)$$
$$\leq{\sf E}_{\bf 0}\exp\{\theta\left(\langle{\bf S}(t),{\bf r}\rangle-t(\nu+\varepsilon)\right)+\alpha_1\beta_1(m_1-1)L(t;{\bf 0})\}=e^{-t\theta(\nu+\varepsilon)}\kappa(t).$$
Here $\theta>0$ and $\kappa(t)={\sf E}_{\bf 0}\exp\{\theta\langle{\bf S}(t),{\bf r}\rangle+\alpha_1\beta_1(m_1-1)L(t;{\bf 0})\}$.
Using properties of conditional expectation we get
$${\sf E}_{\bf 0}\exp\{\theta\langle{\bf S}(t),{\bf r}\rangle+\alpha_1\beta_1(m_1-1)L(t;{\bf 0})\}\mathbb{I}(\tau_{\bf 0}+\overline{\tau}_{{\bf 0},{\bf 0}}\leq t)$$
$$={\sf E}_{\bf 0}\left({\sf E}_{\bf 0}\left(\exp\{\theta\langle{\bf S}(t),{\bf r}\rangle+\alpha_1\beta_1(m_1-1)L(t;{\bf 0})\}|\tau_{\bf 0},\overline{\tau}_{{\bf 0},{\bf 0}}\right)\mathbb{I}(\tau_{\bf 0}+\overline{\tau}_{{\bf 0},{\bf 0}}\leq t)\right)$$
$$=\int_0^t\kappa(t-u)\left(\int_0^u qe^{q(\alpha_1 m_1-1)(u-z)/(1-\alpha_1)}\,d\overline{F}_{{\bf 0},{\bf 0}}(z)\right)\,du=\kappa\ast\eta(t)$$
where sign $\ast$ denotes the convolution of functions and $\eta(u)$
stands for the integral inside big brackets in the previous formula.
We also take into account that $\beta_1=q/(1-\alpha_1)$. Therefore,
\begin{eqnarray*}
\kappa(t)&=&{\sf E}_{\bf 0}\exp\{\theta\langle{\bf S}(t),{\bf r}\rangle+\alpha_1\beta_1(m_1-1)L(t;{\bf 0})\}\mathbb{I}(\tau_{\bf 0}+\overline{\tau}_{{\bf 0},{\bf 0}}>t)\\
&+&{\sf E}_{\bf 0}\exp\{\theta\langle{\bf S}(t),{\bf r}\rangle+\alpha_1\beta_1(m_1-1)L(t;{\bf 0})\}\mathbb{I}(\tau_{\bf 0}+\overline{\tau}_{{\bf 0},{\bf 0}}\leq t)\\
&=&\zeta(t)+\kappa\ast\eta(t)
\end{eqnarray*}
where
\begin{eqnarray*}
\zeta(t)&:=&{\sf E}_{\bf 0}\exp\{\theta\langle{\bf S}(t),{\bf r}\rangle+\alpha_1\beta_1(m_1-1)L(t;{\bf 0})\}\mathbb{I}(\tau_{\bf 0}+\overline{\tau}_{{\bf 0},{\bf 0}}>t)\\
&=&{\sf E}_{\bf 0}\exp\{\theta\langle{\bf S}(t),{\bf r}\rangle+\alpha_1\beta_1(m_1-1)\min\{\tau_{\bf 0},t\}\}\mathbb{I}(\tau_{\bf 0}+\overline{\tau}_{{\bf 0},{\bf 0}}>t).
\end{eqnarray*}

Consider $\theta>1$. Then by convexity of function $H$ the strict inequality $H(\theta{\bf r})>H({\bf r})=\nu$ holds true. Set $\tilde{\kappa}(t)=e^{-tH(\theta{\bf r})}\kappa(t)$, $\tilde{\zeta}(t)=e^{-tH(\theta{\bf r})}\zeta(t)$ and $\tilde{\eta}(t)=e^{-tH(\theta{\bf r})}\eta(t)$. Thus, we get the renewal equation
\begin{equation}\label{renewal_equation}
\tilde{\kappa}(t)=\tilde{\zeta}(t)+\tilde{\kappa}\ast\tilde{\eta}(t).
\end{equation}

By virtue of the definition of the Malthusian parameter one has
$$\alpha_1m_1G^{\ast}_1(\nu)+(1-\alpha_1)G^{\ast}_1(\nu)\overline{F}^{\ast}_{{\bf 0},{\bf 0}}(\nu)=1$$
and, consequently,
$$\overline{F}^{\ast}_{{\bf 0},{\bf 0}}(\nu)=\frac{1-\alpha_1m_1G^{\ast}_1(\nu)}{(1-\alpha_1)G^{\ast}_1(\nu)}=\frac{(1-\alpha_1)\nu-\alpha_1m_1q+q}{q(1-\alpha_1)}.$$
The latter equalities imply that
$$\int_0^{\infty}e^{-\nu u}\eta(u)\,du=\overline{F}^{\ast}_{{\bf 0},{\bf 0}}(\nu)\int_0^{\infty}qe^{-(\nu-q(\alpha_1m_1-1)/(1-\alpha_1))u}\,du=1$$
and, hence,
\begin{equation}\label{eta(t)}
\int_0^{\infty}\tilde{\eta}(u)\,du<\int_0^{\infty}e^{-\nu u}\eta(u)\,du=1.
\end{equation}
In passing we have derived a simple and useful inequality
\begin{equation}\label{simple_useful_inequality}
\nu+q>\alpha_1\beta_1(m_1-1).
\end{equation}

It is not difficult to check with the help of relation \eqref{S(t)=representation} and identity ${\sf E}u^{P(t)}=\exp\{qt(u-1)\}$, $u\in[0,1]$, $t\geq0$, that the stochastic process $\{e^{\theta\langle{\bf S}(t),{\bf r}\rangle-tH(\theta{\bf r})},t\geq0\}$ is a martingale (with respect to filtration ${(\mathcal{F}_t,t\geq0)}$).
In particular,
\begin{equation}\label{e^tH(theta_h(hamma))}
{\sf E}_{\bf 0}e^{\theta\langle{\bf S}(t),{\bf r}\rangle}={\sf E}_{\bf 0}\prod_{i=1}^{P(t)}e^{\theta\langle{\bf Y}^i,{\bf r}\rangle}={\sf E}_{\bf 0}\left(\frac{H(\theta{\bf r})}{q}+1\right)^{P(t)}=e^{tH(\theta{\bf r})}.
\end{equation}
Define measure ${\sf P}^{\theta}$ by martingale change of measure
$$\frac{d{\sf P}^{\theta}}{d{\sf P}_{\bf 0}}=e^{\theta\langle{\bf S}(t),{\bf r}\rangle-tH(\theta{\bf r})}\;\mbox{on}\; \mathcal{F}_t.$$ Then
$$\tilde{\zeta}(t)={\sf E}^{\theta}e^{\alpha_1\beta_1(m_1-1)\min\{\tau_{\bf 0},t\}}\mathbb{I}(\tau_{\bf 0}+\overline{\tau}_{{\bf 0},{\bf 0}}>t)$$
$$={\sf E}^{\theta}e^{\alpha_1\beta_1(m_1-1)t}\mathbb{I}(\tau_{\bf 0}>t)+{\sf E}^{\theta}e^{\alpha_1\beta_1(m_1-1)\tau_{\bf 0}}\mathbb{I}(\tau_{\bf 0}\leq t,\tau_{\bf 0}+\overline{\tau}_{{\bf 0},{\bf 0}}>t)$$
$$=e^{\alpha_1\beta_1(m_1-1)t}{\sf P}^{\theta}(\tau_{\bf 0}>t)+{\sf E}^{\theta}e^{\alpha_1\beta_1(m_1-1)\tau_{\bf 0}}\mathbb{I}(\tau_{\bf 0}\leq t)\left(1-\overline{F}^{\,\theta}_{{\bf 0},{\bf 0}}(t-\tau_{\bf 0})\right)$$
where $\overline{F}^{\,\theta}_{{\bf 0},{\bf 0}}(t):={\sf P}^{\theta}(\overline{\tau}_{{\bf 0},{\bf 0}}\leq t)$, $t\geq0$.
Let us find the distribution of $\tau_{\bf 0}$ with respect to measure ${\sf P}^{\theta}$. Namely, in view of
\eqref{e^tH(theta_h(hamma))} one has
$${\sf P}^{\theta}(\tau_{\bf 0}\leq t)={\sf E}_{\bf 0}\mathbb{I}(\tau_{\bf 0}\leq t)e^{\theta\langle{\bf S}(t),{\bf r}\rangle-tH(\theta{\bf r})}$$
$$=e^{-tH(\theta{\bf r})}{\sf E}_{\bf 0}\left(\mathbb{I}(\tau_{\bf 0}\leq t){\sf E}_{\bf 0}\left(\left.e^{\theta\langle{\bf S}(t),{\bf r}\rangle}\right|\tau_{\bf 0}\right)\right)$$
$$=e^{-tH(\theta{\bf r})}\int_0^t\left({\sf E}_{\bf 0}e^{\theta\langle{\bf S}(t-u),{\bf r}\rangle}\right)\left(\sum_{{\bf x}\in\mathbb{Z}^d,{\bf x}\neq{\bf 0}}e^{\theta\langle{\bf x},{\bf r}\rangle}\frac{q({\bf 0},{\bf x})}{q}\right)qe^{-qu}\,du$$
$$=(H(\theta{\bf r})+q)\int_0^te^{-u(H(\theta{\bf r})+q)}\,du,$$
i.e. $\tau_{\bf 0}$ has an exponential distribution with parameter $H(\theta{\bf r})+q$ with respect to measure
${\sf P}^{\theta}$. Based on this fact we deduce that
$$\tilde{\zeta}(t)=e^{-(H(\theta{\bf r})+q-\alpha_1\beta_1(m_1-1))t}+\left(H(\theta{\bf r})+q\right)\times$$
$$\times\int_0^te^{-(H(\theta{\bf r})+q-\alpha_1\beta_1(m_1-1))u}\left(1-\overline{F}^{\,\theta}_{{\bf 0},{\bf 0}}(t-u)\right)du.$$
Letting $t$ tend to infinity we get
$$\tilde{\zeta}(t)\to\tilde{\zeta}(\infty)=\frac{\left(H(\theta{\bf r})+q\right){\sf P}^{\theta}(\overline{\tau}_{{\bf 0},{\bf 0}}=\infty)}{H(\theta{\bf r})+q-\alpha_1\beta_1(m_1-1)}$$
whenever $H(\theta{\bf r})+q>\alpha_1\beta_1(m_1-1)$. The latter inequality is valid by virtue of \eqref{simple_useful_inequality}. Now to check the estimate $\tilde{\zeta}(\infty)>0$ we have to show that ${{\sf P}^{\theta}(\overline{\tau}_{{\bf 0},{\bf 0}}=\infty)>0}$.

Employing characteristic functions technique one can verify that the process ${({\bf S}(t),t\geq0)}$ has also independent increments with respect to the measure ${\sf P}^{\theta}$.
Moreover, on account of \eqref{e^tH(theta_h(hamma))} one has
\begin{eqnarray*}
{\sf E}^{\theta}S_i(t)&=&{\sf E}_{\bf 0}\left(S_i(t)e^{\theta\langle{\bf S}(t),{\bf r}\rangle-t H(\theta{\bf r})}\right)=\frac{{\sf E}_{\bf 0}\left(S_i(t)e^{\theta\langle{\bf S}(t),{\bf r}\rangle}\right)}{{\sf E}_{\bf 0}e^{\theta\langle{\bf S}(t),{\bf r}\rangle}}\\
&=&\left.\frac{1}{\theta}\frac{\partial\log{\sf E}_{\bf 0}e^{\theta\langle{\bf S}(t),{\bf s}\rangle}}{\partial s_i}\right|_{{\bf s}={\bf r}}=\left.\frac{t}{\theta}\frac{\partial H(\theta{\bf s})}{\partial s_i}\right|_{{\bf s}={\bf r}}.
\end{eqnarray*}
Let us show that ${\sf E}^{\theta}{\bf S}(t)=({\sf E}^{\theta}S_1(t),\ldots,{\sf E}^{\theta}S_d(t))\neq{\bf 0}$, $t>0$. Assume the contrary that $\nabla H(\theta{\bf r})={\bf 0}$. Since the Hessian of $H$ is positive definite, function $H$ reaches the global minimum at point $\theta{\bf r}$. However, ${H(\theta{\bf r})>H({\bf r})=\nu}$. We get the contradiction. Hence,
${\sf E}^{\theta}{\bf S}(t)\neq{\bf 0}$ for each $t>0$, and ${\bf S}$ is a random walk with respect to the measure ${\sf P}^{\theta}$ with non-zero drift. Then the law of large numbers applied to ${\bf S}$ as a process with independent increments entails ${\sf P}^{\theta}(\overline{\tau}_{{\bf 0},{\bf 0}}=\infty)>0$.

Thus, applying the renewal theorem (see, e.g., Theorem~1 in \cite{Feller_71}, Ch.~11, Sec.~6) to renewal equation \eqref{renewal_equation} and taking into account \eqref{eta(t)} we come to relation $$\tilde{\kappa}(t)\to\tilde{\kappa}(\infty)=\frac{\tilde{\zeta}(\infty)}{1-\int_0^{\infty}\tilde{\eta}(u)\,du}\in(0,\infty),\quad t\to\infty.$$

Therefore, if $tH(\theta{\bf r})-t\theta(\nu+\varepsilon)<0$, i.e. $\theta(\nu+\varepsilon)>H(\theta{\bf r})$, then
\begin{eqnarray}
{\sf P}_{\bf 0}\left(\exists v\in N(t):{\bf X}_v(t)\in t\mathcal{O}_{{\bf r},\varepsilon}\right)&\leq&{\sf E}_{\bf 0}\sum\limits_{v\in N(t)}\mathbb{I}\{{\bf X}_v(t)\in t\mathcal{O}_{{\bf r},\varepsilon}\}\label{P(A^c_t)estimate}\\
&\leq&e^{-t(\theta(\nu+\varepsilon)-H(\theta{\bf r}))}\tilde{\kappa}(t).\nonumber
\end{eqnarray}
Denote by $A_{t}$ the event $\{\omega:\forall v\in N(t)\;\mbox{one
has}\;{\bf X}_v(t)\notin t\mathcal{O}_{{\bf r},\varepsilon}\}$. As
usual, $A^c$ stands for the complement of a set $A$ and
$\{A_n\;\mbox{infinitely
often}\;\}=\{A_n\;\mbox{i.o.}\}=\cap_{k=1}^{\infty}\cup_{n=k}^{\infty}A_n$,
for a sequence of sets $A_n$. By virtue of Borel Cantelli's lemma
estimate \eqref{P(A^c_t)estimate} entails ${\sf P}_{\bf
0}\left(A^c_{n/2^m}\;\mbox{i.o.}\right)=0$, for any fixed
$m\in\mathbb{N}$. Consequently, ${\sf P}_{\bf
0}\left(\cap_{m=1}^{\infty}\cup_{k=1}^{\infty}\cap_{n=k}^{\infty}A_{n/2^m}\right)=1$.
It means that for almost all $\omega\in\Omega$ and for any
$m\in\mathbb{N}$ there exists positive integer $k=k(m,\omega)$ such
that for any $n\geq k$ and any $v\in N(n/2^m)$ one has
$X_v(n/2^m)\notin n/2^m\mathcal{O}_{{\bf r},\varepsilon}$. Since the
set of binary rational numbers is dense in $\mathbb{R}$ and the
sojourn time of a particle $v\in N(t)$ in a set $t\mathcal{O}_{{\bf
r},\varepsilon}$ contains non-zero interval with probability $1$, we
conclude that
\begin{equation}\label{P(O_gamma,epsilon)=1}
{\sf P}_{\bf 0}\left(\omega:\exists t_0(\omega)\;\mbox{such that}\;\forall t\geq t_0(\omega)\,\mbox{and}\,\forall v\in N(t),\;{\bf X}_v(t)\notin t\mathcal{O}_{{\bf r},\varepsilon}\right)=1,
\end{equation}
for any $\varepsilon>0$. The assertion \eqref{P(O_gamma,epsilon)=1} remains in force when $\theta$ tends to $1$. Moreover, as $\theta\to 1$ the condition $\theta(\nu+\varepsilon)>H(\theta{\bf r})$ transforms into the trivial one $\nu+\varepsilon>\nu$.

Unfix ${\bf r}\in\mathcal{R}$. If the set $\mathcal{R}$ is finite (it occurs when $d=1$), put $\Upsilon=\mathcal{R}$. Otherwise, let $\Upsilon$ be the everywhere dense set in $\mathcal{R}$ (for instance, let $\Upsilon$ be the set of vectors ${\bf r}$ from $\mathcal{R}$ with rational coordinates $r_1,\ldots,r_{d-1}$). Consider the domain ${\mathcal{O_{\varepsilon}}=\cup_{{\bf r}\in\Upsilon}\mathcal{O_{{\bf r},\varepsilon}}}=\{{\bf x}\in\mathbb{R}^d:{\langle{\bf x},{\bf r}\rangle>\nu+\varepsilon}\;\mbox{for at least one}\;{\bf r}\in\mathcal{R}\}$. Relation \eqref{P(O_gamma,epsilon)=1} entails
$${\sf P}_{\bf 0}\left(\omega:\exists t_1(\omega)\;\mbox{such that}\;\forall t\geq t_1(\omega)\,\mbox{and}\,\forall v\in N(t),\;{\bf X}_v(t)\notin t\mathcal{O}_{\varepsilon}\right)=1.$$
Thus, we obtain the first assertion of Theorem \ref{T:spread} in the case of CBRW with a single catalyst at ${\bf 0}$ and the starting point ${\bf 0}$.

\vskip0.2cm {\it Step 2.} At the second step we also assume that $W=\{{\bf w}_1\}$ with ${\bf w}_1={\bf 0}$ and the starting point of CBRW is ${\bf 0}$. Moreover, we concentrate on the case ${\sf E}\xi^2_1<\infty$. Let us establish Theorem~\ref{T:one_point} and statement~\eqref{T:assertion_2} of Theorem~\ref{T:spread} under these assumptions.

Fix ${\bf r}\in\mathcal{R}$. Let $\delta$ be a number such that $0<\delta<1$. In view of Theorem~4 in \cite{B_Doklady_15} on the set $\mathcal{I}$ at time $\delta t$ there are at least $[Ce^{\nu\delta t}]$ particles at ${\bf 0}$ for some positive constant $C$ (as usual, $[u]$ stands for the integer part of a number $u\in\mathbb{R}_+$).
If these particles move according to the random walk ${\bf S}$ such that $\langle{\bf S}(u),{\bf r}\rangle>0$ for each $u\in[\tau_{\bf 0},t(1-\delta)]$, then far particles in CBRW at time $t$ are not less far than $[Ce^{\nu\delta t}]$ i.i.d. copies of ${\bf S}(t(1-\delta))$ with $\langle{\bf S}(u),{\bf r}\rangle>0$, for each $u\in[\tau_{\bf 0},t(1-\delta)]$. A large deviation estimate (see, e.g., Theorem 4.9.5 of \cite{Borovkov_book_13}) yields, for each $\varepsilon\in(0,\nu)$,
$${\sf P}_{\bf 0}(\langle{\bf S}(u),{\bf r}\rangle>0,u\in[\tau_{\bf 0},t(1-\delta)],\langle{\bf S}(t(1-\delta)),{\bf r}\rangle\geq(\nu-\varepsilon)t)$$
$$=e^{-t(1-\delta)K_{{\bf r},\varepsilon}+o(t)},\quad t\to\infty.$$
Here $K_{{\bf r},\varepsilon}=\inf\left\{\int_0^1L_{\bf r}\left(\varphi'(u)\right)\,du\right\}$ and the infimum is taken over all absolutely continuous functions $\varphi:[0,1]\mapsto\mathbb{R}$ such that $\varphi(0)=0$, $\varphi(u)>0$, ${u\in(0,1)}$, and $\varphi(1)=(\nu-\varepsilon)(1-\delta)^{-1}$. In its turn, function $L_{\bf r}(\theta):=\sup\limits_{\vartheta\in\mathbb{R}}\left(\theta\vartheta-H(\vartheta{\bf r})\right)$, $\theta\in\mathbb{R}$, is the Fenchel-Legendre transform of $H(\vartheta{\bf r})$, $\vartheta\in\mathbb{R}$. The infimum $K_{{\bf r},\varepsilon}=L_{\bf r}\left((\nu-\varepsilon)(1-\delta)^{-1}\right)$ is achieved when $\varphi=\varphi_0$ is a linear function, i.e. $\varphi_0(u)=(\nu-\varepsilon)(1-\delta)^{-1}u$, $u\in[0,1]$, since by Jensen's inequality one has
$$\int_0^1L_{\bf r}\left(\varphi'(u)\right)\,du\geq L_{\bf r}\left(\int_0^1\varphi'(u)\,du\right)
=L_{\bf r}(\varphi(1)-\varphi(0))$$
$$=L_{\bf r}\left((\nu-\varepsilon)(1-\delta)^{-1}\right)=\int_0^1L_{\bf r}\left(\varphi'_0(u)\right)\,du.$$

Letting $\mathcal{Q}_{{\bf r},\varepsilon}:=\{{\bf x}\in\mathbb{R}^d:\langle{\bf x},{\bf r}\rangle<\nu-\varepsilon\}$ (here $0<\varepsilon<\nu$) we get
\begin{equation}\label{P(B^c_t)estimate}
{\sf P}_{\bf 0}\left({\bf X}_v(t)\in t\mathcal{Q}_{{\bf r},\varepsilon}\;\mbox{for any}\;v\in N(t),\,\mu(t;{\bf 0})\geq Ce^{\nu\delta t}\right)
\end{equation}
$$\leq(1-{\sf P}_{\bf 0}(\langle{\bf S}(u),{\bf r}\rangle\!>\!0,u\in[\tau_{\bf 0},t(1\!-\!\delta)],\langle{\bf S}(t(1\!-\!\delta)),{\bf r}\rangle\!\geq\!(\nu\!-\!\varepsilon)t))^{[Ce^{\nu\delta t}]}$$
$$\leq\exp\left\{-\left[Ce^{\nu\delta t}\right] e^{-t(1-\delta)K_{{\bf r},\varepsilon}+o(t)}\right\}=\exp\left\{-e^{(\nu\delta-(1-\delta)K_{{\bf r},\varepsilon})t+o(t)}\right\},\; t\to\infty.$$
Denote by $B_{t}$ the event $\{\omega:\exists v\in N(t)\;\mbox{such
that}\;{\bf X}_v(t)\notin t\mathcal{Q}_{{\bf r},\varepsilon}\}$. By
Borel-Cantelli's lemma estimate \eqref{P(B^c_t)estimate} entails
${\sf P}_{\bf
0}\left(\left.B^c_{n/2^m}\;\mbox{i.o.}\right|\mathcal{I}\right)=0$,
for any fixed $m\in\mathbb{N}$, whenever
\begin{equation}\label{condition_2}
\nu\delta-(1-\delta)L_{\bf r}\left((\nu-\varepsilon)(1-\delta)^{-1}\right)>0.
\end{equation}
Therefore, ${\sf P}_{\bf 0}\left(\left.\cap_{m=1}^{\infty}\cup_{k=1}^{\infty}\cap_{n=k}^{\infty}B_{n/2^m}\right|\mathcal{I}\right)=1$. It means that for almost all $\omega\in\mathcal{I}$ and for any $m\in\mathbb{N}$ there exists positive integer $k=k(m,\omega)$ such that for each $n\geq k$ one can find $v\in N(n/2^m)$ such that $X_v(n/2^m)\notin n/2^m\mathcal{Q}_{{\bf r},\varepsilon}$. Since the set of binary rational numbers is dense in $\mathbb{R}$ and the sojourn time of a particle $v\in N(t)$ in a set $t\mathcal{Q}^c_{{\bf r},\varepsilon}$ contains non-zero interval with probability~$1$, we conclude that
\begin{equation}\label{P(R_gamma,epsilon)=1}
{\sf P}_{\bf 0}\!\left(\left.\omega:\exists t_0(\omega)\;\mbox{such that}\;\forall t\geq t_0(\omega)\,\mbox{one has}\,\exists v\in N(t),\;{\bf X}_v(t)\notin t\mathcal{Q}_{{\bf r},\varepsilon}\right|\mathcal{I}\right)\!=\!1,
\end{equation}
for any $\varepsilon\in(0,\nu)$.

Let us show that, for each $\varepsilon\in(0,\nu)$, there exists $\delta=\delta({\bf r},\varepsilon)\in(0,1)$ such that condition \eqref{condition_2} is satisfied. Indeed, set $a({\bf r})=\left.\frac{\partial H(\theta{\bf r})}{\partial\theta}\right|_{\theta=1}$. Then according to the properties of the Fenchel-Legendre transform (see, e.g., \cite{Borovkov_book_13}, Ch.~1, Sec.~1) we have $L_{\bf r}(a({\bf r}))=a({\bf r})-H({\bf r})=a({\bf r})-\nu\geq0$. It follows that, for $\delta({\bf r},\varepsilon)=1-(\nu-\varepsilon)/a({\bf r})$, inequality \eqref{condition_2} is reduced to the trivial one $\nu>\nu-\varepsilon$. Thus, condition \eqref{condition_2} holds true with $\delta({\bf r},\varepsilon)=1-(\nu-\varepsilon)/a({\bf r})$.

Combination of the proved part of Theorem~\ref{T:spread} and formula \eqref{P(R_gamma,epsilon)=1} implies the assertion of Theorem~2 for the case of a single catalyst at ${\bf 0}$ and the starting point ${\bf 0}$ whenever ${\sf E}\xi^2_1<\infty$. Under the same conditions statement \eqref{T:assertion_2} of Theorem~\ref{T:spread} is established since relation \eqref{P(R_gamma,epsilon)=1} entails
$${\sf P}_{\bf 0}\left(\left.\omega:\exists t_0(\omega)\;\mbox{such that}\;\forall t\geq t_0(\omega)\,\mbox{one has}\,\exists v\in N(t),\;{\bf X}_v(t)\notin t\mathcal{Q}_{\varepsilon}\right|\mathcal{I}\right)=1,$$
for each $\varepsilon\in(0,\nu)$.

\vskip0.2cm {\it Step 3.} At the third step we assume that $W=\{{\bf
w}_1\}$ with ${\bf w}_1={\bf 0}$ and the starting point of CBRW is
${\bf 0}$ whereas now ${\sf E}\xi^2_1=\infty$. To verify assertion
of Theorem~\ref{T:one_point} and statement~\eqref{T:assertion_2} of
Theorem~\ref{T:spread} under such assumptions one can follow the
proof scheme proposed in \cite{Carmona_Hu_2014}, Sec.~5.3, based on
a coupling. It is worthwhile to note that contrast to
\cite{Carmona_Hu_2014} we employ Theorem~3 of \cite{B_Doklady_15}
devoted to the strong convergence of the total and local particles
numbers in supercritical CBRW instead of using properties of a
fundamental martingale as in \cite{Carmona_Hu_2014}. Moreover, here
we exploit function $g(u)=\alpha
f_1\left(q_{esc}+(1-q_{esc})u\right)+(1-\alpha)q_{esc}-u$,
${u\in[0,1]}$, where $q_{esc}={\sf P}_{\bf
0}\left(\overline{\tau}_{{\bf 0},{\bf
0}}=\infty\right)=1-\overline{F}_{{\bf 0},{\bf 0}}(\infty)$ is the
escape probability of the random walk ${\bf S}$. Other details of
the Step~3 proof can be omitted.

\vskip0.2cm {\it Step 4.} Now we consider a supercritical CBRW on $\mathbb{Z}^d$ with a finite catalysts set $W$ and the starting point ${\bf w}_i\in W$. In this case the verification of Theorems~\ref{T:spread} and \ref{T:one_point} repeats mainly the arguments of Steps~1,2 and 3. Therefore we discuss only some differences in these proofs.

Modifying the Step~1 we deal with a system of renewal equations
instead of single renewal equation, namely,
$\kappa_i(t)=\zeta_i(t)+\sum_{j=1}^N\eta_{i,j}\ast\kappa_j(t)$,
$i=1,\ldots,N$, $t\geq0$, where
$$\kappa_i(t)={\sf E}_{{\bf w}_i}\exp\left\{\theta\langle{\bf S}(t),{\bf r}\rangle+\sum\nolimits_{j=1}^N\alpha_j\beta_j(m_j-1)L(t;{\bf w}_j)\right\},$$
\begin{eqnarray*}
\zeta_i(t)&=&{\sf E}_{{\bf w}_i}\exp\left\{\theta\langle{\bf S}(t),{\bf r}\rangle+\alpha_i\beta_i(m_i-1)\min\left\{\tau_{{\bf w}_i},t\right\}\right\}\\
&\times&\mathbb{I}(\tau_{{\bf w}_i}+\,_{W_j}\overline{\tau}_{{\bf w}_i,{\bf w}_j}>t,\,j=1,\ldots,N),
\end{eqnarray*}
$$\eta_{i,j}(t)=\int_0^t qe^{q(\alpha_i m_i-1)(t-u)/(1-\alpha_i)}\,d\,{_{W_j}\overline{F}_{{\bf w}_i,{\bf w}_j}(u)}.$$

Denoting by $J(\lambda)$ and $K(\lambda)$ matrices with the corresponding entries $\int_0^{\infty}{e^{-\lambda u}\eta_{i,j}(u)\,du}$ and $\delta_{i,j}\left(1-\alpha_im_iq/\left(\lambda(1-\alpha_i)+q\right)\right)$, $i,j=1,\ldots,N$, $\lambda\geq\nu$, one can check the following identity, for each $\rho\in\mathbb{R}$,
$$D(\lambda)-\rho I=K(\lambda)(J(\lambda)-\rho I)$$
where, as usual, $I$ is the identity matrix.
Hence, since for $\lambda\geq\nu$ diagonal matrix $K(\lambda)$ is non-degenerate, irreducible matrix $J(\lambda)$ has the Perron root $\rho(J(\lambda))$ (a positive eigenvalue of maximal modulus with respect to other eigenvalues of the matrix) equal to $1$ if and only if $\lambda=\nu$. It follows that the Perron root of matrix $J\left(H(\theta{\bf r})\right)$ (when $\theta>1$) is strictly less than $1$.

In the same manner, as in Step~1, one can derive that
$$\tilde{\zeta}_i(t)\to\tilde{\zeta}_i(\infty)=\frac{\left(H(\theta{\bf r})+q\right)\left(1-\sum_{j=1}^N {_{W_j}\overline{F}^{\,\theta}_{{\bf w}_i,{\bf w}_j}(\infty)}\right)}{H(\theta{\bf r})+q-\alpha_i\beta_i(m_i-1)},\quad t\to\infty,$$
and the finite limit $\left(\tilde{\zeta}_1(\infty),\ldots,\tilde{\zeta}_N(\infty)\right)$ is not identically zero. Then applying the renewal theorem (see, e.g., Theorem~2.2, item (ii), of \cite{Mode_68_2}) to the system of renewal equations $\tilde{\kappa}_i(t)=\tilde{\zeta}_i(t)+\sum_{j=1}^N\tilde{\eta}_{i,j}\ast\tilde{\kappa}_j(t)$, $i=1,\ldots,N$, $t\geq0$, we come to relation $\tilde{\kappa}_i(t)\to\tilde{\kappa}_i(\infty)>0$, for each $i=1,\ldots,N$, as $t\to\infty$, with
$$\left(\tilde{\kappa}_1(\infty),\ldots,\tilde{\kappa}_N(\infty)\right)=\left(\tilde{\zeta}_1(\infty),\ldots,\tilde{\zeta}_N(\infty)\right)\left(I-J\left(H(\theta{\bf r})\right)^{\top}\right)^{-1}$$
where $ ^{\top}$ means the matrix transposition. The rest of the proof of statement~\eqref{T:assertion_1} in case of CBRW with general catalysts set $W$ and the starting point from $W$ as well as the verification of statement~\eqref{T:assertion_2} and Theorem~\ref{T:one_point} is implemented similar to arguments of Steps~1,2 and 3.

\vskip0.2cm {\it Step 5.} Turning to a supercritical CBRW on
$\mathbb{Z}^d$ with a finite catalysts set $W$ and the starting
point ${\bf x}\notin W$, we supplement the catalysts set $W$ with
${\bf w}_{N+1}={\bf x}$ and put $\alpha_{N+1}=0$, $m_{N+1}=0$,
$G_{N+1}(t)=1-e^{-qt}$, $t\geq0$. According to Lemma~3 of
\cite{B_TPA_15} a new CBRW with catalysts set $\{{\bf
w}_1,\ldots,{\bf w}_{N+1}\}$ is supercritical whenever the
underlying CBRW is supercritical and the Malthusian parameters in
these CBRW coincide. Then one can apply the proved parts of
Theorems~\ref{T:spread} and \ref{T:one_point} to the new CBRW and
obtain the desired assertions of those theorems for CBRW with an
arbitrary starting point.

\vskip0.2cm Thus, the proof of Theorems~\ref{T:spread} and \ref{T:one_point} is complete.

\section{Proof of Theorem~\ref{T:front} and Examples}\label{s:front} Firstly, we prove Theorem~\ref{T:front}. To this end put
$\mathcal{Z}=\{{\bf z}({\bf r}):{\bf r}\in\mathcal{R}\}$. Let us
verify inclusion $\mathcal{Z}\subset\mathcal{P}$. Indeed, according
to Theorem~23.5 of \cite{Rockafellar_70}, for any ${\bf r},{\bf
r}'\in\mathcal{R}$ one has
\begin{equation}\label{L:relation}
\left\langle{\bf z}({\bf r}),{\bf r}'\right\rangle=\nu\frac{\left\langle\nabla H({\bf r}),{\bf r}'\right\rangle}{\left\langle\nabla H({\bf r}),{\bf r}\right\rangle}\leq\nu\frac{L(\nabla H({\bf r}))+H\left({\bf r}'\right)}{L(\nabla H({\bf r}))+H({\bf r})}=\nu
\end{equation}
where function $L({\bf s})$, ${\bf s}\in\mathbb{R}^d$, is the Fenchel-Legendre transform of function $H({\bf s})$, ${\bf s}\in\mathbb{R}^d$, and the sign $\leq$ transforms into the sign $=$ if and only if ${\bf r}'={\bf r}$. In other words, we established that $\left\langle{\bf z}({\bf r}),{\bf r}\right\rangle=\nu$ and $\left\langle{\bf z}({\bf r}),{\bf r}'\right\rangle<\nu$, ${\bf r}'\neq{\bf r}$. Thus, $\mathcal{Z}\subset\mathcal{P}$.

Now we check that $\mathcal{P}\subset\mathcal{Z}$. Let ${\bf x}\in\mathcal{P}$. Since the spherical image of a closed convex surface is a unit sphere (see, e.g., \cite{Busemann_1958}, Ch.~1, Sec.~4) and the function $H$ is smooth, one can find $\sigma>0$ and ${\bf r}\in\mathcal{R}$ such that ${\bf x}=\sigma\nabla H({\bf r})$. Therefore, by virtue of \eqref{L:relation} one has
$\left\langle{\bf x},{\bf r}'\right\rangle=\sigma\left\langle\nabla H({\bf r}),{\bf r}'\right\rangle\leq\sigma\langle\nabla H({\bf r}),{\bf r}\rangle$
where the sign $\leq$ transforms into the sign $=$ if and only if ${\bf r}'={\bf r}$. Since ${\bf x}\in\mathcal{P}$, we deduce that $\sigma=\nu\langle\nabla H({\bf r}),{\bf r}\rangle^{-1}$, i.e. ${\bf x}={\bf z}({\bf r})\in\mathcal{Z}$.

Thus, the proof of Theorem~\ref{T:front} is complete.

\vskip0.2cm Now consider five examples.

\vskip0.2cm {\it Example 1.} Focus on a continuous-time counterpart of the discrete-time CBRW on $\mathbb{Z}$ treated in \cite{Carmona_Hu_2014}. Then the set $\mathcal{R}$ consists of two points $r_1$ and $r_2$ being the roots of equation $H(r)=\nu$. Since $H(0)=0$ and $H$ is a convex function, we see that $r_1<0<r_2$. Hence, the propagation front $\mathcal{P}$ also consists of two points $\nu\,r_1^{-1}$ and $\nu\,r_2^{-1}$, i.e. for large time $t$ all the particles are located almost surely at set $t\left(\mathcal{P}\cup\mathcal{Q}\right)=t\{x\in\mathbb{R}:\nu r_1^{-1}\leq x\leq\nu r_2^{-1}\}$. This conclusion implies a result of \cite{Carmona_Hu_2014}.

\begin{figure}
\includegraphics[width=17cm]{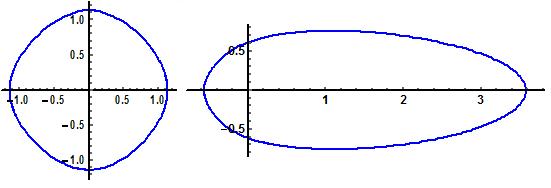}
\caption[]{To the left and to the right the corresponding plots of $\mathcal{P}$ for Examples 2a and 2b.}\label{Plots_example2ab}
\end{figure}

\vskip0.2cm {\it Example 2a.} Concentrate on the simplest case of CBRW on $\mathbb{Z}^2$, i.e. when jumps of the random walk occur to the neighboring points with probabilities $1/4$. Then $H({\bf s})=q(e^{s_1}+e^{s_2}+e^{-s_1}+e^{-s_2})/4-q=q(\cosh{s_1}+\cosh{s_2})/2-q$, ${\bf s}\in\mathbb{R}^2$.
Solving equation $H(r_1,r_2)=\nu$ with respect to unknown variable $r_2$ we obtain $r_2=\pm\arch\left(2\nu q^{-1}+2-\cosh r_1\right)$ where ${r_1\in\left[-\arch\left(2\nu q^{-1}+1\right),\arch\left(2\nu q^{-1}+1\right)\right]}$. Consequently,
$$\nabla H({\bf r})=\left(q\sinh r_1/2,\pm q\sinh\!\left(\arch\!\left(2\nu q^{-1}+2-\cosh r_1\right)\!\right)\!/2\right)\;\mbox{and}\;{\bf z}({\bf r})$$
$$\!=\!\!\frac{\nu\left(\sinh r_1,\pm\sinh\left(\arch\left(2\nu q^{-1}+2-\cosh r_1\right)\right)\right)}{r_1{\sinh r_1}\!+\!\arch\!\left(2\nu q^{-1}+2\!-\!\cosh r_1\right)\!\sinh\!\left(\arch\!\left(2\nu q^{-1}+2\!-\!\cosh r_1\right)\!\right)},$$
for ${\bf r}\in\mathcal{R}$. The plot of $\mathcal{P}$ is drawn on Figure \ref{Plots_example2ab} to the left when $\nu=2$ and $q=2$.

\vskip0.2cm {\it Example 2b.} Consider now non-symmetric CBRW on $\mathbb{Z}^2$, i.e., for instance, when the random walk instantly shifts to the vector $(2,0)$, $(-1,0)$, $(0,1)$ and $(0,-1)$ with corresponding probabilities $1/2$, $1/6$, $1/6$ and $1/6$. Then $H({\bf s})=q\left(e^{2s_1}/2+e^{-s_1}/6+\cosh s_2/3\right)-q$, ${\bf s}\in\mathbb{R}^2$. Solving equation ${H(r_1,r_2)=\nu}$ with respect to $r_2$ we get $r_2=\pm\arch\left(3\nu q^{-1}+3-3e^{2r_1}/2-e^{-r_1}/2\right)$ where $r_1\in\mathbb{R}$ is such that $3e^{3r_1}-2\left(3\nu q^{-1}+2\right)e^{r_1}+1\leq0$. It follows that
$\nabla H({\bf r})=\left(q\left(e^{2r_1}-e^{-r_1}/6\right),\pm q\sinh\left(\arch\left(3\nu q^{-1}+3-3e^{2r_1}/2-e^{-r_1}/2\right)\right)/3\right)$. In a similar way one can write exact formula for ${\bf z}({\bf r})$, ${\bf r}\in\mathcal{R}$. The plot of $\mathcal{P}$ is represented on Figure \ref{Plots_example2ab} to the right when $\nu=1$ and $q=3$.

\begin{figure}
\includegraphics[width=17cm]{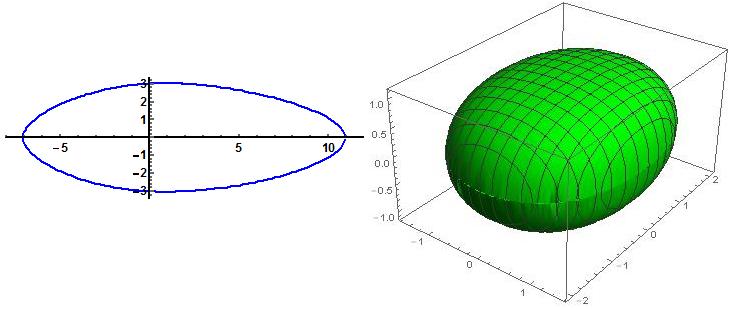}
\caption[]{To the left and to the right the corresponding plots of $\mathcal{P}$ for Examples 2c and 3.}\label{Plots_example2c3}
\end{figure}

\vskip0.2cm {\it Example 2c.} Let us concentrate on a non-symmetric CBRW on $\mathbb{Z}^2$ with non-bounded jump sizes. Namely, let random walk instantly shift to the vector $(n,0)$, $(-n,0)$, $(0,1)$ and $(0,-1)$, $n\in\mathbb{N}$, with corresponding probabilities $\sigma_1^{n-1}e^{-\sigma_1}/{(4(n-1)!)}$, $\sigma_2^{n-1}e^{-\sigma_2}/{(4(n-1)!)}$, $1/4$ and $1/4$. In other words, the jumps of the random walk in the right and in the left directions obey the displaced Poisson law with parameters $\sigma_1>0$ and $\sigma_2>0$, respectively, whereas the shifts of the random walk to the top or to the bottom occur to neighboring points only. Then
\begin{eqnarray*}
H({\bf s})&=&q\left(\frac{e^{-\sigma_1}}{4}\sum_{n=1}^{\infty}\frac{e^{s_1n}\sigma_1^{n-1}}{(n-1)!}
+\frac{e^{-\sigma_2}}{4}\sum_{n=1}^{\infty}\frac{e^{-s_1n}\sigma_2^{n-1}}{(n-1)!}+\frac{\cosh s_2}{2}\right)-q\\
&=&q\left(\frac{e^{\sigma_1\left(e^{s_1}-1\right)+s_1}}{4}+\frac{e^{\sigma_2\left(e^{-s_1}-1\right)-s_1}}{4}+\frac{\cosh s_2}{2}\right)-q,\quad{\bf s}\in\mathbb{R}^2.
\end{eqnarray*}
Solving equation ${H(r_1,r_2)=\nu}$ with respect to unknown variable $r_2$ we come to equality $r_2=\pm\arch\left(2\nu q^{-1}+2-e^{\sigma_1\left(e^{r_1}-1\right)+r_1}/2-e^{\sigma_2\left(e^{-r_1}-1\right)-r_1}/2\right)$ where ${r_1\in\mathbb{R}}$ is such that the argument of function $\arch$ in the previous formula is not less than $1$. Hence,
\begin{eqnarray*}
&&\nabla H({\bf r})=\left(\frac{q}{4}\left(e^{\sigma_1\left(e^{r_1}-1\right)+r_1}\left(\sigma_1 e^{r_1}+1\right)-e^{\sigma_2\left(e^{-r_1}-1\right)-r_1}\left(\sigma_2 e^{-r_1}+1\right)\right),\right.\\
&&\left.\pm\frac{q}{2}\sinh\left(\arch\left(2\nu q^{-1}+2-\frac{1}{2}e^{\sigma_1\left(e^{r_1}-1\right)+r_1}-\frac{1}{2}e^{\sigma_2\left(e^{-r_1}-1\right)-r_1}\right)\right)\right)
\end{eqnarray*}
and the precise formula for ${\bf z}({\bf r})$, ${\bf r}\in\mathcal{R}$, can be written in a similar way. The plot of $\mathcal{P}$ is drawn on Figure \ref{Plots_example2c3} to the left when $\nu=4$, $q=8$, $\sigma_1=2$ and $\sigma_2=1$.

\vskip0.2cm {\it Example 3.} Consider CBRW on $\mathbb{Z}^3$ such that coordinates of the random walk are independent and its jump ${\bf Y}=(Y_1,Y_2,Y_3)$ has the following marginal distributions. For $n\in\mathbb{N}$, set
$${\sf P}(Y_1=n)={\sf P}(Y_1=-n)=\frac{\sigma_1^{n-1}e^{-\sigma_1}}{2(n-1)!},$$
$${\sf P}(Y_2=n)={\sf P}(Y_2=-n)=\frac{\sigma_2^{n-1}e^{-\sigma_2}}{2(n-1)!},\quad{\sf P}(Y_3=1)={\sf P}(Y_3=-1)=\frac{1}{2}.$$
Since function $H({\bf s})$ can be represented in the form $H({\bf s})=q\left({\sf E}e^{\left\langle{\bf s},{\bf Y}\right\rangle}-1\right)$, one has
$$H({\bf s})=q\frac{V(\sigma_1,s_1)}{2}\cdot\frac{V(\sigma_2,s_2)}{2}\cdot\cosh s_3-q,\quad{\bf s}\in\mathbb{R}^3,$$
where $V(\sigma,u):=e^{\sigma\left(e^{u}-1\right)+u}+e^{\sigma\left(e^{-u}-1\right)-u}$,
$\sigma>0$ and $u\in\mathbb{R}$. Solving equation ${H(r_1,r_2,r_3)=\nu}$ with respect to unknown variable $r_3$ we get
\begin{equation}\label{example3_r3}
r_3=\pm\arch\frac{4\nu q^{-1}+4}{V(\sigma_1,r_1)V(\sigma_2,r_2)}
\end{equation}
where $r_1,r_2\in\mathbb{R}$ are such that the argument of function $\arch$ in formula \eqref{example3_r3} is not less than $1$.
Then
$$\nabla H({\bf r})=(\nu+q)\left(\frac{e^{\sigma_1\left(e^{r_1}-1\right)+r_1}\left(\sigma_1 e^{r_1}+1\right)-e^{\sigma_1\left(e^{-r_1}-1\right)-r_1}\left(\sigma_1 e^{-r_1}+1\right)}{V(\sigma_1,r_1)}\right.,$$
$$\frac{e^{\sigma_2\left(\!e^{r_2}\!-\!1\right)\!+r_2}\!\left(\sigma_2e^{r_2}\!+\!1\right)\!-\!e^{\sigma_2\left(\!e^{-r_2}\!-\!1\right)\!-\!r_2}\!\left(\sigma_2
e^{\!-\!r_2}\!+\!1\right)\!}{V(\sigma_2,r_2)},\left.\frac{V(\sigma_1,r_1)V(\sigma_2,r_2)}{4\nu q^{-1}+4}\sinh r_3\right)$$
where $r_3$ is described by relation \eqref{example3_r3}. Similarly one can write the explicit formula for ${\bf z}({\bf r})$, ${\bf r}\in\mathcal{R}$. The plot of $\mathcal{P}$ is represented on Figure \ref{Plots_example2c3} to the right when $\nu=0.5$, $q=1$, $\sigma_1=0.2$ and $\sigma_2=0.5$.


\begin{thebibliography}{99}

\bibitem{AB_00} Albeverio S. and Bogachev L.V. Branching random walk in a catalytic medium.
I. Basic equations. \emph{Positivity} {\bf 4}(2000), no.~1, 41-100.

\bibitem{Bertacchi_Zucca_15} Bertacchi D. and Zucca F. Branching random walks and multi-type contact-processes on the percolation cluster of $\mathbb{Z}^d$. \emph{Ann. Appl. Probab.} {\bf 25}(2015), no.~4, 1993-2012.

\bibitem{Biggins_78} Biggins J.D. The asymptotic shape of the branching random walk. \emph{Adv. Appl. Probab.} {\bf 10}(1978), no.~1, 62-84.

\bibitem{Bocharov_Harris_14} Bocharov S. and Harris S.C. Branching Brownian motion with catalytic branching at the origin. \emph{Acta Appl. Math.} {\bf 134}(2014), no.~1, 201-228.

\bibitem{Borovkov_book_13} Borovkov A.A. \emph{Asymptotic Analysis of Random Walks}. Moscow, FIZMATLIT, 2013 (in Russian).

\bibitem{Bremaud_99} Br\'{e}maud P. \emph{Markov chains: Gibbs Fields, Monte-Carlo Simulation, and Queues}. Springer, New York, 1999.

\bibitem{B_TrudyMIAN_13} Bulinskaya E.Vl. Subcritical catalytic branching random walk with
finite or infinite variance of offspring number. \emph{Proc. Steklov
Inst. Math.} {\bf 282}(2013), no.~1, 62-72.

\bibitem{B_JTP_14} Bulinskaya E.Vl. Local particles numbers in critical
branching random walk. \emph{J. Theoret. Probab.} {\bf 27}(2014), no.~3, 878-898.

\bibitem{B_SPL_14} Bulinskaya E.Vl. Finiteness of hitting times under
taboo. \emph{Statist. Probab. Lett.} {\bf 85}(2014), no.~1, 15-19.

\bibitem{B_TPA_15} Bulinskaya E.Vl. Complete classification of catalytic branching processes. \emph{Theory Probab. Appl.} {\bf 59}(2015), no.~4, 545-566.

\bibitem{B_Doklady_15} Bulinskaya E.Vl. Strong and weak convergence of the population size in a supercritical catalytic branching process. \emph{Doklady Math.} {\bf 92}(2015), no.~3, 714-718.

\bibitem{Busemann_1958} Busemann H. \emph{Convex Surfaces}. Interscience Publishers, New York, 1958.

\bibitem{Comets_Popov_07} Comets F. and Popov S. On multidimensional branching random walks in random environment. \emph{Ann. Probab.} {\bf 35}(2007), no.~1, 68-114.

\bibitem{Carmona_Hu_2014} Carmona Ph. and Hu Y. The spread of a catalytic branching random
walk. \emph{Ann. Inst. Henri Poincar\'{e} Probab. Stat.} {\bf 50}(2014), no.~2, 327-351.

\bibitem{CKMV_09} Cranston M., Korallov L., Molchanov S. and Vainberg B. Continuous model for homopolymers. \emph{J. Funct. Anal.} {\bf 256}(2009), no.~8, 2656-2696.

\bibitem{DR_13} Doering L. and Roberts M. Catalytic branching processes via
spine techniques and renewal theory. In: Donati-Martin C., et al.
(Eds.), \emph{S\'{e}minaire de Probabilit\'{e}s} XLV, \emph{Lecture
Notes in Math.} {\bf 2078}(2013), 305-322.

\bibitem{Feller_71} Feller W. \emph{An Introduction to Probability Theory and Its Applications}. Vol.II. Wiley, New York, 1971.

\bibitem{GKM_07} G\"{a}rtner J., Konig W. and Molchanov S. Geometric characterization of intermittency in the parabolic Anderson model. \emph{Ann. Probab.}, {\bf 35}(2007), no.~2, 439-499.

\bibitem{Harris_Roberts_16} Harris S. and Roberts M. The many-to-few lemma and multiple spines.
\emph{Ann. Inst. Henri Poincar\'{e} Probab. Stat.} (to appear), available at http://arxiv.org/abs/1106.4761.

\bibitem{HJV_05} Haccou P., Jagers P. and Vatutin V.A. \emph{Branching Processes: Variation, Growth, and Extinction of Populations}. Cambridge University Press, Cambridge, 2005.

\bibitem{HTV_12} Hu Y., Topchii V.A. and Vatutin V.A. Branching random
walk in ${\bf Z}^{4}$ with branching at the origin only.
\emph{Theory Probab. Appl.} {\bf 56}(2012), no.~2, 193-212.

\bibitem{KA_15} Kimmel M. and Axelrod D. \emph{Branching Processes in Biology}. Springer, New York, 2015.

\bibitem{Lawler_Limic_10} Lawler G. and Limic V. \emph{Random Walk: a Modern Introduction}. Cambridge University Press, Cambridge, 2010.

\bibitem{Le-Gall_Lin_15} Le Gall J.-F. and Lin S. The range of tree-indexed random walk in low dimensions. \emph{Ann. Probab.} {\bf 43}(2015), no.~5, 2701-2728.

\bibitem{Mallein_15} Mallein B. Maximal displacement of $d$-dimensional branching Brownian motion. \emph{Electron. Commune. Probab.} {\bf 20}(2015), no.~76, 1-12.

\bibitem{Mode_68_2} Mode Ch.J. A multidimensional age-dependent branching process with applications to natural selection. II. \emph{Math. Biosci.} {\bf 3}(1968), 231-247.

\bibitem{Molchanov_Yarovaya_2012} Molchanov S.A. and Yarovaya E.B. Branching processes with lattice spatial dynamics and a finite set of particle generation centers. \emph{Doklady Math.} {\bf 446}(2012), no.~3, 259-262.

\bibitem{Rockafellar_70} Rockafellar R.T. \emph{Convex Analysis}. Princeton University Press, Princeton, 1970.

\bibitem{Shi_LNM_15} Shi Z. \emph{Branching Random Walks}. \'{E}cole d'\'{E}t\'{e} de Probabilit\'{e}s de Saint-Flour XLII - 2012, \emph{Lecture Notes in Math.} {\bf 2151}(2015).

\bibitem{VTY} Vatutin V.A., Topchii V.A. and Yarovaya E.B. Catalytic branching random walk and
queueing systems with random number of independent servers.
\emph{Theory Probab. Math. Statist.} (2004), no.~69, 1-15.

\bibitem{Y_TPA_11} Yarovaya E.B. Criteria of exponential growth
for the numbers of particles in models of branching random walks.
\emph{Theory Probab. Appl.} {\bf 55}(2011), no.~4, 661-682.

\end{thebibliography}
\end{document}